\documentclass[12pt, reqno]{amsart}

\usepackage{fullpage}
\usepackage{graphicx} 
\usepackage{amsthm,amssymb,amsmath}
\usepackage{mathrsfs}
\usepackage{tikz}
\usepackage{euscript}
\usepackage{wasysym}
\usepackage[margin=1in]{geometry}
\usepackage{bbm}

\usepackage[hypertexnames=false]{hyperref}
\hypersetup{colorlinks=true,
	citecolor=blue,
	filecolor=blue,
	linkcolor=blue,
}

\usepackage[hypertexnames=false]{hyperref}
\hypersetup{colorlinks=true,
	citecolor=blue,
	filecolor=blue,
	linkcolor=blue,
}

\title{Deranged Perfect Matchings on complete graph and balanced complete r-partite graph}

\author{Boqing Deng}

\newtheorem{theorem}{Theorem}[section]

\newtheorem{definition}[theorem]{Definition}

\newtheorem{lemma}[theorem]{Lemma}

\newtheorem{corollary}[theorem]{Corollary}

\newcommand\pma{\ensuremath{\mathrm{pm}}}
\newcommand\bpm{\ensuremath{\mathrm{bpm}}}
\newcommand\po{\ensuremath{\mathrm{Po}}}

\newcommand{\pr}{\ensuremath{\mathbb P}}
\newcommand{\vn}[1]{\left|\boldsymbol{\mathbf #1}\right|_1}

\newcommand{\bd}[1]{\boldsymbol{\mathbf #1}}

\newcommand{\dseq}{(D _m)_{m=1}^\ell}

\newcommand{\dtv}{d_{TV}(\bd X, \bd Y)}

\newcommand{\cg}{\ensuremath{K_{2n}}}
\newcommand{\cbg}{\ensuremath{K_{r \times (r-1)n}}}

\newcommand{\xseq}{(X_m)_{m=1}^\ell}
\newcommand{\yseq}{(Y_m)_{m=1}^\ell}
\newcommand{\psa}{\ensuremath{\psi_1}}
\newcommand{\psb}{\ensuremath{\psi_2}}
\newcommand{\psc}{\ensuremath{\psi_3}}
\newcommand{\pha}{\ensuremath{\varphi_1}}
\newcommand{\phb}{\ensuremath{\varphi_2}}
\makeatletter
\newcommand*\bigcdot{\mathpalette\bigcdot@{.5}}
\newcommand*\bigcdot@[2]{\mathbin{\vcenter{\hbox{\scalebox{#2}{$\m@th#1\bullet$}}}}}
\makeatother

\begin{document}

\begin{abstract}
We proved that for any finite collection of sparse subgraphs $\dseq$ of the complete graph $\cg$, and a uniformly chosen perfect matching $R$ in $\cg$, the random vector $(|E(R \cap D_m)|)_{m=1}^\ell$ jointly converges to a vector of independent Poisson random variables with mean $|E(D_m)|/(2n)$. We also showed a similar result when $\cg$ is replaced by the balanced complete $r$-partite graph $K_{r \times 2n/r}$ for fixed $r$ and determined the asymptotic joint distribution. The proofs rely on elementary tools of the Principle of Inclusion-Exclusion and generating functions. These results extend recent works of Johnston, Kayll and Palmer, Spiro and Surya, and Granet and Joos from the univariate to the multivariate setting.
\end{abstract}
\maketitle
\section{Introduction}
For a sequence simple graphs $G_n$ where $\lim_{n \to \infty} V(G_n) = \infty$, let $\pma(G_n)$ denote the number of perfect matchings of $G_n$. If $M_n$ is an arbitrary perfect matching of $G_n$, the problem of determining the asymptotic ratio
\[
\frac{\pma(G_n-M_n)}{\pma(G_n)}
\]
when $n \to \infty$ is of great combinatorial interest. For example, when $G_n = K_{n,n}$ is the balanced complete bipartite graphs, the limiting ratio equals the number of permutations with no fixed point and converges to the limit $e^{-1}$. A permutation without a fixed point is called a derangement, and hence a perfect matching of $G_n-M_n$ is called a deranged perfect matching. Counting the number of derangement was proposed by Montmort \cite{Montmort1708-13} in 1708, and this problem can be solved using the principle of Inclusion-Exclusion (See \cite{Stanley} for more details). When $G_n$ is the complete graph $K_{2n}$, the limiting ratio equals $e^{-1/2}$. Brawner \cite{Brawner} conjectured this asymptotic ratio, and it was proved by Margolius \cite{Margolius}. For the rest of this paper, we will omit the index $n$ and write $G_n, M_n$ as $G,M$. 

Let $r = r(n)$ be a integer valued function and $r| 2n$, and let $K_{r \times 2n/r}$ denote the balanced complete $r$-partite graph, where $V(K_{r \times 2n/r}) = V_1 \uplus ... \uplus V_r$ is a partition of the vertex set,  $|V_i| = 2n/r$ for each $i \in [r]$, and there is an edge between $u$ and $v$ if and only if $u, v$ lies in different parts. As a generalization of both $K_{n,n}$ and $K_{2n/r}$, Johnston, Kayll and Palmer \cite{JKP} conjectured that when $G = K_{r \times 2n/r}$, the limit ratio converges to $e^{-r/(2r-2)}$. They solved the conjecture when $r$ is linearly proportional to $n$, $r = \Omega(n^{\delta})$ for some $\delta > 0$, and a simplified variant when $r$ is constant.

Spiro and Surya \cite{SS} fully solved the previous conjecture and proved that when $R$ is a uniformly random perfect matching of $K_{r \times 2n/r}$, the number of edges in $R \cap M$ converges to the Poisson distribution with mean $r/(2r-2)$. 

Granet and Joos \cite{GJ} generalized $G$ to regular robust expander graphs and $M$ to any matchings or spanning regular subgraphs. Suppose $R$ is a uniformly random perfect matching of $G$, and $D$ is a matching or a regular spanning subgraph of $G$. In that case, they showed the number of edges $R$ intersecting with $D$ converges to the Poisson distribution with mean $|E(D)|/\deg(G)$. 

In this paper, we generalize the result of \cite{JKP, SS, GJ} by extending the distribution into a multivariate joint Poisson distribution on the $\ell$ dimensional integer lattice for some fixed $\ell$. We generalized the graph $D$ to a collection of subgraphs $\dseq$, and we proved that the asymptotic joint distribution of $(|E(R \cap D_m)|)_{m=1}^\ell$ is a multivariate Poisson distribution. Specificaly, if $\dseq$ are disjoint, then the asymptotic distribution of $(|E(R \cap D_m)|)_{m=1}^\ell$ is independent. This phenomenon provides macroscopic evidence supporting the heuristic proposed by Granet and Joos, which we will outline in Section \ref{Sec2}. We shall define the distance of total variation of two random vectors, taking values in the $\ell$ dimensional integer lattice as follows.
\begin{definition}
Let $\bd X, \bd Y$ be two random vectors taking values in $\mathbb N^\ell$. We denote the distance of total variation of $\bd X$ and $\bd Y$ by
\[
\dtv = \sum_{\bd k \in \mathbb N^\ell}\left|\pr(\bd X = \bd k) - \pr(\bd Y = \bd k)\right|
\]
\end{definition}
Section \ref{Sec3} treats the complete graph $K_{2n}$ as the parent graph and also examines $K_{2n}-N$, where $N$ is a sparse subgraph. A specific example of $K_{2n}-N$ is when $G = K_{r \times 2n/r}$ where $r|2n$ and $r$ is linearly proportional to $n$. In this case, $\deg(G) = n - n/r$ where $n/r$ is a constant independent of $n$. The next theorem states the main result in Section \ref{Sec3}.
\begin{theorem}
    \label{MainThm1}
    Let $\ell, C$ be fixed constants independent of $n$, $N$ be a subgraph of $\cg$ and $\dseq$ be a collection of disjoint subgraphs of $\cg - N$ such that $\Delta(N), \Delta(D_m) \leq C$ for all $m$ and some constant $C$. Let $R$ be a uniformly random perfect matching of $\cg-N$. Let $\bd X = \xseq, \bd Y = \yseq$ be two random vectors such that $X_m = |E(R \cap D_m)|, Y_m$ indepdendently follows $\po(|E(D_m)|/2n)$. Then,
    \[
    \lim_{n \to \infty} \dtv = 0
    \]
\end{theorem}

As mentioned before, Johnston, Kayll and Palmer \cite{JKP} proposed and solved a simplified varient of determining the asymptotic ratio 
$\pma(K_{r \times 2n/r} - M)/\pma(K_{r \times 2n/r})$. They defined a balanced perfect matching of $K_{r \times 2n/r}$ as a perfect matching such that the number of edges between $V_i$ and $V_j$ is the same for all $i \neq j \in [r]$. Let $\bpm(\cdot)$ denote the function that counts the number of balanced perfect matchings. If $r$ is constant on $n$ and $M$ is a balanced perfect matching of $K_{r \times 2n/r}$, Johnston, Kayll and Palmer \cite{JKP} proved the following result:
\begin{equation}
\label{JKPeqt1}
 \lim_{n \to \infty} \frac{\bpm(K_{r \times 2n/r}-M)}{\bpm(K_{r \times 2n/r})} = e^{-r/(2r-2)}   
\end{equation}

In Section \ref{Sec4}, we generalize (\ref{JKPeqt1}) and derive an analogous theorem to Theorem \ref{MainThm1}. To ensure the existence of a balanced perfect matching in a balanced complete $r$-partite graph, we must have $r(r-1) | 2n$. Therefore, for convenience we replace $2n$ by $r(r-1)n$. The next Theorem is the main result of Section \ref{Sec4}.
\begin{theorem}
    \label{MainThm2}
    Let $\ell, C, r$ be fixed, $R$ be a uniformly random \textbf{balanced} perfect matching of $\cbg$, $\dseq$ be a collection of disjoint subgraphs of $\cbg$ such that $\Delta(D_m) \leq C$ for all $m$ and for some constant $C$. Let $\bd X = \xseq, \bd Y = \yseq$ be two random vectors such that $X_m = |E(R \cap D_m)|, Y_m$ independently follows $\po(E(D_m)/(r-1)^2n)$. Then,
    \[
    \lim_{n \to \infty} \dtv = 0
    \]
\end{theorem}

\subsection{Conventions and Notations}
\label{Sec1.1}
Throughout the discussion in the paper, we assume that the dimensions $\ell, \lambda$ and the constants $r, C$ are fixed, independent of $n$. Objects such as the graphs $G, D_m$ and the random vectors $\bd X, \bd Y$, are defined with respect to $n$. We say a quantity is fixed or constant if it is independent of $n$ and unless stated otherwise, a quantity depends on $n$.

We denote the natural number $\mathbb N = \{0,1,2,...,\}$ and for each $k \geq 1 \in \mathbb N$, we denote $[k] = \{1,2,...,k-1,k\}$. For $a, b \in \mathbb N$, we define the falling factorial $a_{(b)} = a(a-1)...a(a-b+1)$. We will use the standard Landau notations $o(\cdot),O(\cdot)$ etc. We assume these asymptotic notations are defined with respect to $n$ when $n \to \infty$.

We define a generating function as an analytic power series from $\mathbb C^\ell$ to $\mathbb C$. Specifically, we define a generating function $G: \mathbb C^{\ell} \to \mathbb C$ as a power series
\[
G(\bd s) = \sum_{\bd k \in \mathbb N^\ell} \alpha_{\bd x} \bd s^{\bd k}
\]
that converges absolutely for all $\bd k \in \mathbb C^{\ell}$. 

For simplicity, in the rest of the paper we will use the multi-index notation to simplify the presentation. The object to consider are $\ell \times \lambda$ dimensional indices $\bd x = (x_{m,k})_{m \in [\ell], k \in [\lambda]} \in \mathbb N^{\ell \times \lambda}$ and $\ell$ dimensional variables $\bd s = (s_m)_{m=1}^\ell \in \mathbb N^\ell$. For any vector $\bd \alpha = (\alpha_1, ..., \alpha_d), \bd \beta = (\beta_1, ..., \beta_d) \in \mathbb N^d, \gamma \in \mathbb N$ of dimension $d$, we define
\begin{align*}
  \bd \alpha + \bd \beta &= (\alpha_1 + \beta_1, ..., \alpha_d + \beta_d) \qquad \text{and} \ \bd \alpha - \bd \beta \ \text{similarly}\\
  \gamma \bd \alpha &= (\gamma \alpha_1, ..., \gamma \alpha_d) \\
  \bd \alpha^{\bd \beta} & = \alpha_1^{\beta_1}...\alpha_d^{\beta_d} \\
  \bd \alpha! &= \alpha_1!...\alpha_d! \\
  \binom{\bd \alpha}{\bd \beta} & = \binom{\alpha_1}{\beta_1}...\binom{\alpha_d}{\beta_d} = \frac{\bd \alpha!}{\bd \beta! (\bd \alpha-\bd \beta)!}\\
  \binom{\vn  \alpha }{\bd \alpha} &= \binom{\vn \alpha}{\alpha_1, ..., \alpha_d} =  \frac{\vn \alpha!}{\alpha_1!...\alpha_d!} = \frac{\vn \alpha !}{\bd \alpha!} \qquad \textrm{this is the multinomial coefficient}\\
  \bd \alpha_{(\bd \beta)} &= {\alpha_1}_{(\beta_1)}...{\alpha_d}_{(\beta_d)}
\end{align*}
We denote $\bd 1 = (1,1,...,1)$ to be the all-one-vector where the dimension will be clear from the context. We also define $O(\gamma)\bd 1$ to be the set of $\bd \beta = (\beta_1, ..., \beta_d)$ where $\beta_1, ..., \beta_d = O(\gamma)$. For $\bd x \in \mathbb N^{\ell \times \lambda}$, we define \begin{align*}
   \vn x &= \sum_{m \in [\ell], k \in \lambda} x_{m,k} \in \mathbb N \\
   \bd x_m &= (x_{m,k})_{k=1}^\lambda \in \mathbb N^{\lambda} \\
   |\bd x_m|_1 &= \sum_{k =1}^\lambda x_{m,k} \in \mathbb N
\end{align*}
We define the function $\psa: \mathbb N^{\ell \times \lambda} \to \mathbb N^\ell$ by $\psa(\bd x) = (\vn{x_m})_{m=1}^\ell $

\subsection{Organization}
The remainder of this paper will be organized as follows. In Section \ref{Sec2}, we will give a heuristic reasoning of Theorem \ref{MainThm1} and Theorem \ref{MainThm2}, a sketch of the proof, and a list of tools used in the proof. In Section \ref{Sec3}, we will prove Theorem \ref{MainThm1}. In Section \ref{Sec4}, we will prove Theorem \ref{MainThm2}. In Section \ref{Sec5}, we will generalize Theorem \ref{MainThm1} and Theorem \ref{MainThm2} by dropping the restriction that $\dseq$ are disjoint. In Section \ref{Sec6}, we will suggest potential directions for future works on this problem.

\section{Sketch and Tools of the proof}
\label{Sec2}

\subsection{Sketch of the proof}
For simplicity, if $G = \cg$ is a complete graph and $R$ is a uniformly random perfect matching of $G$, then for any edge $e \in E(G)$, we have $\pr(e \in R) = 1/(2n-1) = (1+o(1))/2n$. If $G = \cbg$, then $\pr(e \in R) = 1/(r-1)^2n$.  If we have $k$ edges, $e_1, ..., e_k$ in $\cg$ or $\cbg$, $R$ is chosen uniformly from $\cg$ or $\cbg$, and $n$ is much larger than $k$, for all $i \in [k]$, the probability that $R$ contains all of $e_i$ for  when $e_i$ are non incident equals.
\[
\pr(\forall i \in [k],\  e_i \in R) = (1+o(1))\prod_{i=1}^k \pr(e_i \in R)
\]
The event that $R$ contains each $e_i$ is roughly independent. Suppose there exists a pair $e_i$, $e_j$ of incident edges, then we have $\pr(\forall i \in [k],\  e_i \in R) = 0$. However, if we uniformly select $k$ edges from a sparse subgraph of $\cg$ or $\cbg$, it is rare that there exists an incidient pair.

Granet and Joos \cite{GJ} suggests a heuristic that if $G$ is a $d$-regular graph and $D$ is a regular sparse subgraph of $G$, then the probability that each edge in $D$ intersects $R$ is roughly independent and identical, and $|E(R \cap D)|$ will approximately follow the binomial distribution $\textrm{Binom}(|E(D)|, 1/d)$, which will converge to $\po(|E(D)|/d)$ as $n \to \infty$. They proved this convergence given that $G$ is a robust expander graph. Particularly, if $G= K_{2n}$, then $d = 2n-1 = (1+o(1))2n$ and if $G = \cbg$, then $d = (r-1)^2n$. If we have a collection of disjoint sparse graphs $D_m$, we expect that the joint distribution $|E(R \cap D_m)|$ should converge to the independent Poisson distribution. 

Our proof strategy is as follows. We introduce the probability generating function $G(s_1, ... , s_{\ell})$ for which the coefficient for $s_1^{r_1}...s_{\ell}^{r_\ell}$ is the probability that $R$ intersect $D_m$ in exactly $r_m$ edges for all $m \in [\ell]$. We generalized the method of Johnston, Kayll, and Palmer \cite{JKP}, by using the principle of Inclusion-Exclusion to estimate the coefficient of this generating function. We will show that the coefficient in this probability generating functions gets close to the coefficient in the probability generating function of independent multivariate Poisson distribution as $n$ gets large, thus showing the two random vectors converge in distance of total variation.

If $G = K_{2n}$, since $X_m$ is roughly independent, we expect the conditional distribution $X_m|X_1 = 0$ to roughly equal the distribution of $X_m$. Therefore, given that $R$ does not intersect $D_1$, the joint distribution $(|E(R \cap D_m)|)_{m=2}^\ell$ should still be independently Poisson. Specifically, the distribution of $(|E(R^* \cap D_m)|)_{m=2}^\ell$ when $R^*$ is chosen uniformly from $\cg - D_1$ should be the same as the distribution of $(|E(R \cap D_m)|)_{m=2}^\ell$ when $R$ chosen uniformly in $\cg$, conditioned on the event that $R$ does not intersect $D_1$.  We will use the notation $N$ to denote this specific sparse subgraph $D_1$, and we hence generalize the base graph $G$ from $\cg $ to $ \cg - N$. When $r$ is constant and the parent graph is $G = \cbg$, we will use a more sophisticated counting argument to determine the probability generating function, but the core idea and the structure of the proof are similar.

Finally, when $D_m$ is not disjoint, the joint distribution $X_m$ is no longer independent. We can break $D_m$ into disjoint parts, where each part's intersection to $R$ follows the Poisson distribution. Since a sum of independent random variables, each with a Poisson distribution still has a Poisson distribution, we can determine the limiting joint distribution as a not necessarily independent Poisson joint distribution. We will describe the process of decomposing graphs in Section \ref{Sec5}. 
\subsection{Tools for the proof}
The two main theorems we use are the Principle of Inclusion-Exclusion and Tannery's Theorem. Both are also used in the work of Johnston, Kayll, and Palmer \cite{JKP}. We present the Principle of Inclusion-Exclusion in the form of generating functions. Interested readers could refer to \cite{GF} for reference. The following results presents the principle of Inclusion-Exclusion in the fullest generality needed for the proof.

 The next theorem is a version of the Principle of Inclusion-Exclusion. For our application, $\mathcal U$ will be the set of perfect matchings of the parent graph $G$, $I_{m,k}$ will be a partition of the edges of each graph in $\dseq$ in $k$ parts, and the $m,k$\ th coordinate of $\bd P$ will be the set of edges of intersection between the perfect matching and $I_{m,k}$.

 To better articulate the next theorem and arguments in the rest of the paper, we develop some notations here. For a class of sets $\bd S = (S_{m,k})_{m \in [\ell], k \in [\lambda]}$, we define function $\pha$ such that
\[
\pha(\bd S) = (\sum_k |S_{m,k}|)_{m=1}^\ell \in \mathbb N^{\ell}
\]

For two class of sets $\bd S = (S_{m,k})_{m \in [\ell], k \in [\lambda]}, \bd T = (T_{m,k})_{m \in [\ell], k \in [\lambda]}$, we say $\bd S \sqsubset \bd T$ if $S_{m,k} \subset T_{m,k}$ for each $m \in [\ell],k \in [\lambda]$. Let the notation $\mathcal P(\cdot)$ denote the power set.
\begin{theorem}[Principle of Inclusion-Exclusion]
\label{pie}

Let $\mathcal U$ be a finite universal set, $\bd I = (I_{m,k})_{m \in [\ell], k \in [\lambda]}$ be a finite collection of finite index sets. Let $\bd P: \mathcal U \to \prod_{m,k} \mathcal P(I_{m,k})$ be a function. For each $\bd S \sqsubset \bd I$, we define the number
\[
N(\sqsupset \bd S) = |\{\omega \in \mathcal U: \bd S \sqsubset \bd P(\omega)\}|
\]
Then, for each $\bd r \in \mathbb N^\ell$, the coefficient of $\bd s^{\bd r}$ in the generating function
\[
G(\bd s) = \sum_{\bd S}N(\sqsupset \bd S)(\bd s - \bd 1)^{\pha(\bd S)}
\]
is the number of $\omega$ such that $\pha(\bd P(\omega)) = \bd r$.    
\end{theorem}
\begin{proof}
We define the function 
\[
\phb(\bd S) = (|S_{m,k}|)_{m \in [\ell], k \in [\lambda]} \in \mathbb N^{\ell \times \lambda}
\]
We know
\begin{align*}
 G(\bd s + \bd 1) &= \sum_{\bd S}N(\sqsupset \bd S)\bd s^{\pha(\bd S)} \\
&= \sum_{\bd S}\sum_{\omega: \bd S \sqsubset \bd P(\omega)} \bd s^{\pha(\bd S)} \\
&= \sum_{\omega}\sum_{\bd S: \bd S \sqsubset \bd P(\omega)}\bd s^{\pha(\bd S)} \\
&= \sum_{\omega}\sum_{\bd x \in \mathbb N^{\ell \times \lambda}} \binom{\phb(\bd P(\omega))}{\bd x}\bd s^{\psa (\bd x)} \\
&= \sum_{\omega}(\bd s + \bd 1)^{\pha(\bd P(\omega))}
\end{align*}
The last equality is due to a version of the Binomial Theorem for vectors in multiple dimensions. It can be proved by applying the binomial theorem separately for each index $m,k$ and multipling the result all together.

The coefficient of $(\bd s + \bd 1)^{\bd r}$ is the number of $\omega$ such that $\pha(\bd P(\omega)) = \bd r$. Therefore, we obtain the original generating function by substituting $\bd s + \bd 1$ with $\bd s$.
\end{proof}
Theorem \ref{pie} has a straightforward extention in terms of probability-generating function. We shall state this as the next corollary and use this form of the Principle of Inclusion-Exclusion in Section \ref{Sec3} and Section \ref{Sec4}. We shall formally define the probability generating function.
\begin{definition}
    \label{probgf}
    Let $\bd X$ be a random vector on $\mathbb N^\ell$. A probability generating function of $\bd X$ is a generating function defined as
    \[
    G_{\bd X}(\bd s) = \sum_{\bd k} \pr(\bd X = \bd k)\bd s^{\bd k}
    \]
\end{definition}
\begin{corollary}
    \label{pieCor}
    Let $\Omega$ be a finite sample space where each sample is assigned a uniform probability measure. Let $\bd I = (I_{m,k})_{m \in [\ell], k \in [\lambda]}$ be a finite collection of finite index sets. Let $\bd P: \Omega \to \prod_{m,k} \mathcal P(I_{m,k})$ be a function. For each $\bd S \sqsubset \bd I$, define the event $A_{\bd S} = \{\omega \in \Omega: \bd S \sqsubset \bd P(\omega)\}$.  Let $\bd X$ be the $\ell$ dimensional random vector such that $\bd X(\omega) = \pha(\bd P(\omega))$. Then, the probability generating function of $\bd X$ is given by
    \[
    G_{\bd X}(\bd s) = \sum_{\bd S} \pr(A_{\bd S})(\bd s - \bd 1)^{\pha(\bd S)}
    \]
\end{corollary}
For our application, $\Omega$ will be the sample space of all perfect matchings of $G$, and the event $A_{\bd S}$ is the set of perfect matchings $\omega$ of $G$ such that $\omega \cap I_{m,k}$ contains the $m,k$\ th coordinate of $\bd S$ for each $m,k$.

The next theorem is a special case of Lebesgue dominated convergence theorem and provides a sufficient condition for interchanging limits and infinite summation. Interested readers could refer to \cite{analysis} for reference and proof using only elementary mathematical analysis. \\
Let $I$ be a countably infinite set and $\{s_i\}_{i \in I}$ be a sequence of real numbers. If there exists a bijection $g: \mathbb N \to I$ such that $\sum_n s_{g(n)}$ converges absolutely, then we can define the sum $\sum_{i \in I} s_i$ as $\sum_n s_{g(n)}$. In this case, the choice of $g$ will not alter the sum.

\begin{theorem}[Tannery's Theorem]
\label{tannery}
Let $I$ be a countably infinite index set, $\{f_i(n)\}_{i \in I}$ be a sequence of functions from $\mathbb N$ to $\mathbb R$, $\{\alpha_i\}_{i \in I}$, $\{M_i\}_{i \in I}$ be two sequences of real numbers. If for each $i \in I, n \in \mathbb N, |f_i(n)| \leq M_i$, for each $i \in I$, $\lim_{n \to \infty} f_i(n) = \alpha_i$, and if $\sum_{i \in I} M_i < \infty$, then $\sum_{i \in I} f_i(n)$ is defined for each $n$ and we have
\[
\lim_{n \to \infty} \sum_{i \in I} f_i(n) = \sum_{i \in I} \alpha_i
\]
\end{theorem}

For the remainder of this paper, we always choose the set $I$ to be $\mathbb N^{\ell \times \lambda}$, the finite dimensional lattice on natural numbers. The next lemma gives an upper bound of the total variation distance of two random vectors in terms of their probability-generating function.

\begin{lemma}
\label{bounddtv}
Let $\bd X, \bd Y$ be two random vectors on $\mathbb N^{\ell}$. Suppose there exists 
\[
\{\alpha_{\bd x}\}_{\bd x \in \mathbb N^{\ell \times \lambda}}, \{\beta_{\bd x}\}_{\bd x \in \mathbb N^{\ell \times \lambda}} \subset \mathbb R
\]
such that the probability generating function of $\bd X, \bd Y$ satisfies
\begin{align*}
    G_{\bd X}(\bd s) &= \sum_{\bd x \in \mathbb N^{\ell \times \lambda}} \alpha_{\bd x}(\bd s - \bd 1)^{\psa(\bd x)}\\
    G_{\bd Y}(\bd s) &= \sum_{\bd x \in \mathbb N^{\ell \times \lambda}} \beta_{\bd x}(\bd s - \bd 1)^{\psa(\bd x)} \\
\end{align*}
then their total variation distance satisfies
\[
\dtv \leq \sum_{\bd x \in \mathbb N^{\ell \times \lambda}} |\alpha_{\bd x}-\beta_{\bd x}|2^{\vn x}
\]
\end{lemma}
\begin{proof}
We can write
\begin{align*}
  G_{\bd X}(\bd s) &= \sum_{\bd x \in \mathbb N^{\ell \times \lambda}}  \alpha_{\bd x}(\bd s - \bd 1)^{\psa(\bd x)} \\
  &=\sum_{\bd x \in \mathbb N^{\ell \times \lambda}} \alpha_{\bd x} \sum_{\bd k \in \mathbb N^{\ell \times \lambda}} \binom{\bd x}{\bd k}(-\bd 1)^{\psa(\bd x - \bd k)}\bd s^{\psa(\bd k)}\\
  &= \sum_{\bd k \in \mathbb N^{\ell \times \lambda}}\left( \sum_{\bd x \in \mathbb N^{\ell \times \lambda}} \alpha_{\bd x}\binom{\bd x}{\bd k}(-\bd 1)^{\psa(\bd x - \bd k)}\right)\bd s^{\psa(\bd k)}
\end{align*}
Where the second equality is due to the multi-dimensional Binomial Theorem.
Similarly, we can write
\[
G_{\bd Y}(\bd s) = \sum_{\bd k \in \mathbb N^{\ell \times \lambda}}\left( \sum_{\bd x \in \mathbb N^{\ell \times \lambda}} \beta_{\bd x}\binom{\bd x}{\bd k}(-\bd 1)^{\psa(\bd x - \bd k)}\right)\bd s^{\psa(\bd k)}
\]
Hence, by Triangle Inequality
\begin{align*}
    \dtv &= \sum_{\bd k \in \mathbb N^{\ell \times \lambda}}\left|\left( \sum_{\bd x \in \mathbb N^{\ell \times \lambda}} \alpha_{\bd x}\binom{\bd x}{\bd k}(-\bd 1)^{\psa(\bd x - \bd k)}\right) - \left( \sum_{\bd x \in \mathbb N^{\ell \times \lambda}} \beta_{\bd x}\binom{\bd x}{\bd k}(-\bd 1)^{\psa(\bd x - \bd k)}\right)\right| \\
    & \leq \sum_{\bd k \in \mathbb N^{\ell \times \lambda}}\sum_{\bd x \in \mathbb N^{\ell \times \lambda}} |\alpha_{\bd x}-\beta_{\bd x}|\binom{\bd x}{\bd k} \\
    & = \sum_{\bd x \in \mathbb N^{\ell \times \lambda}} |\alpha_{\bd x}-\beta_{\bd x}|2^{\vn x}
\end{align*}
where the last equality is due to the Binomial Theorem.
\end{proof}
The next technical lemma gives a bound of a difference of two products when the difference of each coordinates are small.
\begin{lemma}
    \label{boundprod}
    Let $\{a_i\}_{i \in [d]}, \{b_i\}_{i \in [d]}, K$ be real numbers such that $0 \leq b_i \leq a_i \leq K$ for each $i$. Then
    \[
    \prod_1^d a_i - \prod_1^d (a_i-b_i) \leq K^{d-1}\sum_i b_i
    \]
\end{lemma}
\begin{proof}
    We have
    \begin{align*}
        \prod_1^d a_i - \prod_1^d (a_i-b_i) &= \sum_{j=1}^d \left(\prod_{i<j}(a_i-b_i)\right)\left(\prod_{i\geq j}a_i\right)-\left(\prod_{i\leq j}(a_i-b_i)\right)\left(\prod_{i>j}a_i\right) \\
        &=\sum_{j=1}^d b_i\left(\prod_{i<j}(a_i-b_i)\right)\left(\prod_{i>j} a_i \right) \\
        & \leq \sum_i b_i K^{d-1}
    \end{align*}
    The second inequality is due to the fact that both $a_i$ and $a_i-b_i$ are bounded by $K$, and so the lemma follows.
\end{proof}

\section{Case of $K_{2n}$}
\label{Sec3}
In this section, we assume $\lambda = 1$, $\ell$ to be finite and independent of $n$. Therefore, we denote $\bd x = (x_m)_{m=1}^\ell$, and we denote $\vn x = \sum_{m =1}^\ell x_m$. We begin from the case where $G = \cg$ and $D_m$ are disjoint. The next theorem is the main result in this section, and we will prove it by proving a series of lemmas.
\begin{theorem}
\label{MainThm3.1}
Let $\dseq$ be a collection of disjoint subgraphs of $K_{2n}$, where $\Delta(D_m) \leq C$ for each $m$. Let $R$ be a uniformly random perfect matching of $K_{2n}$. Define $\bd X = (X_m), \bd Y = (Y_m), X_m = |E(R \cap D_m)|$, $Y_m \sim \po(|E(D_m)|/2n)$ independently. Then,
$\lim_{n \to \infty} d_{TV}(\bd X, \bd Y) = 0$.
\end{theorem}
We shall use the following notation to help with counting.
\begin{definition}
    \label{xmatching1}
    Let $\bd x \in \mathbb N^\ell$. We define a $\bd x$-matching of $\cg$ as a $\vn x$-matching $M$ where $|E(M \cap D_m)| = x_m$ for each $m$. We define $\mu_{\bd x}$ as the number of $\bd x$-matchings on $\cg$.
\end{definition}
We prove the following lemma by applying the principle of Inclusion-Exclusion. In this case, we apply Corollary \ref{pieCor}. Since $\lambda = 1$, for a set $\bd S = (S_m)_{m=1}^\ell$, we define $\varphi_1(\bd S) = (|S_m|)_{m=1}^\ell$ for convenience.
\begin{lemma}
    \label{GF3}
    Let $\bd X$ be defined as Theorem \ref{MainThm3.1}. The probability-generating function of $\bd X$ is given by
    \[
    G_{\bd X}(\bd s) = \sum_{\bd x: \vn x \leq n} \mu_{\bd x} \frac{(n)_{(\vn x)}2^{\vn x}}{(2n)_{(2\vn x)}}(\bd s - \bd 1)^{\bd x}
    \]
\end{lemma}
\begin{proof}
    Let $\Omega$ be the set of all perfect matchings of $\cg$ and let $\bd I = (|E(D_m)|)_{m=1}^\ell$. For $\bd S \sqsubset \bd I$, let $\hat{\bd S} = \bigcup_{m=1}^\ell S_m$ be the natural identification of $\bd S$ into a sets of edges in $\cg$. If $\hat{\bd S}$ is a matching, i.e. the edges in $\hat{\bd S}$ are non-incident, it is a $\bd x$-matching where $\bd x = \pha(\bd S)$. Then we must have $\vn x \leq n$ and
    \[
    \pr(A_{\bd S}) = \frac{\pma(K_{2n - 2 \vn x})}{\pma(\cg)} = \frac{(2n-2\vn x)!}{(n-\vn x)!2^{n-\vn x}} \Bigg / \frac{(2n)!}{n!2^n} = \frac{(n)_{\vn x}2^{\vn x}}{(2n)_{(2\vn x)}}
    \]
    Where $(2n-2\vn x)!/(n-\vn x)!2^{n-\vn x}$ is the number of ways to extend $\hat{\bd S}$ into a perfect matching of $\cg$. If $\hat{\bd S}$ is not a matching, then $\pr(A_{\bd S}) = 0$.
    Since there are $\mu_{\bd x}$ $\bd x$-matchings, we know
    \begin{align*}
        \sum_{\bd S}\pr(A_{\bd S})(\bd s - \bd 1)^{\pha(\bd S)} & = \sum_{\bd x : \vn x \leq n} \sum_{\bd S : \varphi(\bd S) = \bd x} \pr(A_{\bd S})(\bd s - \bd 1)^{\bd x} \\
        &= \sum_{\bd x: \vn x \leq n} \mu_{\bd x} \frac{(n)_{\vn x}2^{\vn x}}{(2n)_{(2\vn x)}} (\bd s - \bd 1)^{\bd x}
    \end{align*}
    and the lemma follows.
\end{proof}
For convenience, we write $|E(D_m)| = d_m$ and $\bd d = (d_m)_{m=1}^\ell$. The next lemma estimates the quantity $\mu_{\bd x}$. 
\begin{lemma}
\label{boundmu3}
Let $\mu_{\bd x}, C$ be defined as above. Then,
\[
\mu_{\bd x} = \frac{(\bd d - O(\vn x)\bd 1)^{\bd x}}{\bd x!}
\]
\end{lemma}
\begin{proof}
    We choose a $\bd x$-matching $M_{\bd x}$ one edge at a time. For each $1 \leq m \leq \ell, 1 \leq k \leq  x_m$, suppose that we have chosen the edges in $M_{\bd x} \cap D_{m'}$ for each $m' < m$ and we have chosen $k-1$ edges in $D_m$. Then, we have chosen at most $\vn x$-edges. The chosen edges are adjacent to at most $2C\vn x$ edges, so we can choose the $k$th edge in at least $d_m- 2C\vn x$ ways and at most $d_m$ ways. Applying the product rule for counting and dividing $x_m!$ for each $m$ as we choose the edges unorderly, we know there are
    $(d_m - O(\vn x))/x_m!$ ways to choose the $x_m$ edges in $D_m$ to form the $\bd x$-matching. Applying the product rule again,  we have
    \[
    \mu_{\bd x} = \prod_m \frac{(d_m - O(\vn x))^{x_m}}{x_m!} = \frac{(\bd d - O(\vn x) \bd 1)^{\bd x}}{\bd x!}
    \]
\end{proof}
Recall that $\bd Y = (Y_m)_{m=1}^\ell$ and $Y_m$ follows independently to $\po(d_m/2n)$. Let $\bigcdot : \mathbb R^\ell \times \mathbb R^\ell \to \mathbb R$ denote the vector dot product. We can write the probability generating function of $\bd Y$ as
\begin{align*}
    G_{\bd Y}(\bd s) &= e^{(1/2n)\bd d \cdot (\bd s - 1)} \\
    &= \sum_{\bd x \in \mathbb N^\ell} \frac{(\bd d/2n)^{\bd x}}{\bd x!}(\bd s - 1)^{\bd x}
\end{align*}
Our last step is to apply Tannery's theorem. Set
\[
\alpha_{\bd x} = \begin{cases}
    \mu_{\bd x}\frac{(n)_{(\vn x)}2^{\vn x}}{(2n)_{(2 \vn x)}} &  \vn x \leq n \\
    0 & \text{otherwise}
\end{cases} \qquad \beta_{\bd x} = \left(\prod_m \frac{(|E(D_m)|/2n)^{x_m}}{x_m!} \right)
\]
Then, by Theorem \ref{bounddtv}, we know
\[
d_{TV}(\bd X, \bd Y) \leq \sum_{\bd x \in \mathbb N^\ell}|\alpha_{\bd x} - \beta_{\bd x}|2^{\vn x}
\]
We want to show
\[
\lim_{n \to \infty }\sum_{\bd x \in \mathbb N^\ell}|\alpha_{\bd x} - \beta_{\bd x}|2^{\vn x} = 0
\]
In order to apply Tannery's theorem to switch the order of limit and infinite summation, we want to show that there exists $\gamma_{\bd x}$, independent of $n$, such that 
\[
|\alpha_{\bd x}-\beta_{\bd x}| = |\alpha_{\bd x}(n)-\beta_{\bd x}(n)| \leq \gamma_{\bd x} \qquad \text{for all} \  n \in \mathbb N
\]
and $\sum_{\bd x} \gamma_{\bd x}2^{\vn x} < \infty$. Then, we want to show, $|\alpha_{\bd x} - \beta_{\bd x}| \to 0$ for all $\bd x$ fixed when $n \to \infty$. If so, by Tannery's theorem, we have
\[
\lim_{n \to \infty }\sum_{\bd x \in \mathbb N^\ell}|\alpha_{\bd x} - \beta_{\bd x}|2^{\vn x} = \sum_{\bd x \in \mathbb N^\ell}\lim_{n \to \infty}|\alpha_{\bd x} - \beta_{\bd x}|2^{\vn x} = 0
\]
This proves $\dtv \to 0$.
\begin{proof}[Proof of Theorem \ref{MainThm3.1}]
We will first find $\gamma_{\bd x}$ such that $\sum_{\bd x}\gamma_{\bd x}2^{\vn x} < \infty$. If $\vn x \leq n$, we have
\begin{align}
  \frac{(n)_{(\vn x)}2^{\vn x}}{(2n)_{(2 \vn x)}} &= \frac{2^{\vn x}}{n^{\vn x}} \frac{(n)_{(\vn x)}n^{\vn x}}{(2n)_{(2 \vn x)}} \nonumber \\
  &= \frac{2^{\vn x}}{n^{\vn x}} \left(\prod_{k=1}^{\vn x} \frac{n}{2n-k+1}\right)\left(\prod_{k=1}^{\vn x}\frac{n-k+1}{2n-\vn x - k+1}\right) \nonumber \\
  &\leq \frac{2^{\vn x}}{n^{\vn x}} \label{term3}
\end{align}
as $n \leq 2n-k+1$ and $n-k+1 \leq 2n - \vn x - k + 1$. By Lemma \ref{boundmu3}, $\mu_{\bd x} \leq \frac{\bd d^{\bd x}}{\bd x!}$, therefore, we have
\[
0 \leq \alpha_{\bd x}, \beta_{\bd x} \leq \frac{(2 \bd d/n)^{\bd x}}{\bd x!} \leq \frac{(2C)^{\vn x}}{\bd x!}
\]
hence $\gamma_{\bd x} = (2C)^{\vn x}/\bd x!$ will suffice, and $\gamma_{\bd x}$ is independent of $n$. We have
\[
\sum_{\bd x} \frac{(2C)^{\vn x}}{\bd x!} 2^{\vn x}  = e^{4 \ell C} < \infty
\]
and we have completed the first part. Then, we want to show that for all fixed $\bd x$, $\lim_{n \to \infty}|\alpha_{\bd x}-\beta_{\bd x}| = 0$. We know given $n$ sufficiently large, holding $\bd x$ fixed,
\[
|\alpha_{\bd x}-\beta_{\bd x}| \leq \left|\frac{\bd d^{\bd x}}{\bd x!}- \frac{(\bd d - O(\vn x)\bd 1)^{\bd x}}{\bd x!}\right|\frac{(n)_{(\vn x)}2^{\vn x}}{(2n)_{(2 \vn x)}} + \frac{\bd d^{\bd x}}{\bd x!}\left|\frac{(n)_{(\vn x)}2^{\vn x}}{(2n)_{(2 \vn x)}} - (2n)^{- \vn x}\right|
\]
by triangle inequality. 
We can view $\bd d^{\bd x} - (\bd d - O(\vn x)\bd 1)^{\bd x}$ as a product of $\vn x$ terms since $d_m \leq \vn d \leq Cn$, we can apply Lemma \ref{boundprod} and (\ref{term3}) to obtain the following bound
\begin{equation}
\label{term1}
\left|\frac{\bd d^{\bd x}}{\bd x!}- \frac{(\bd d - O(\vn x)\bd 1)^{\bd x}}{\bd x!}\right|\frac{(n)_{(\vn x)}2^{\vn x}}{(2n)_{(2 \vn x)}} \leq \vn x O(\vn x)(Cn)^{\vn x - 1}\frac{2^{\vn x}}{n^{\vn x}}
\end{equation}
Moreover,
\begin{equation}
 \label{term2}
\frac{\bd d^{\bd x}}{\bd x!}\left|\frac{(n)_{(\vn x)}2^{\vn x}}{(2n)_{(2 \vn x)}} - (2n)^{- \vn x}\right| \leq (Cn)^{\vn x}(2n)^{-\vn x}\frac{(n)_{(\vn x)}2^{\vn x}(2n)^{\vn x} - (2n)_{(2\vn x)}}{(2n)_{(2\vn x)}}   
\end{equation}

Holding $\vn x$ fixed, the numerator of the fraction in (\ref{term2}) is a polynomial of $n$ of degree at most $2\vn x-1$, so the fraction can be bounded by $o(1)$.
Therefore, we can view both (\ref{term1}) and (\ref{term2}) as bounded by $o(1)$, and we have $\lim_{n \to \infty} |\alpha_{\bd x}-\beta_{\bd x}| = 0$. Hence, we have
\[
\lim_{n \to \infty} \dtv = 0
\]
\end{proof}
If we let $\ell = 1$, we obtain the next corollary. It can also be shown using a result of Godsil \cite{Godsil} and Zaslavsky \cite{Zaslavsky}.
\begin{corollary}
    \label{Cor3.2}
    Let $D$ be a subgraph of $\cg$ such that $\Delta(D) \leq C$ for some constant $C$, and $R$ be a uniformly random perfect matching of $\cg$. Then
    \[
    \lim_{n \to \infty} \pr(|E(R \cap D)| = 0) = \lim_{n \to \infty} e^{-|E(D)|/2n}
    \]
\end{corollary}
Note that this ratio is bounded below by $e^{-(C/2)(1+o(1))}$ and is asymptotically nonzero.
In particular, if $D$ is $d$-regular for some constant $d$, then the limiting probability equals $e^{-d/2}$. \\
Now we want to generalize $\cg$ to $\cg - N$ for some graphs $N$ where $\Delta(N) \leq C$. We can define $X^* = |E(R \cap N)|$. Applying previous results, we show that the random vector $\bd X, X^*$ jointly converges to some independent Poisson distribution. Therefore, we expect the conditional distribution $\bd X|X^*$ to be close to $\bd X$ asymptotically. In particular, the distribution $\bd X|X^* = 0$ should also be independently Poisson. The next result is a corollary of Theorem \ref{MainThm3.1} and is a restatement of Theorem \ref{MainThm1}.
\begin{corollary}
    Let $N, (D_m)_{m=1}^\ell$ be disjoint subgraphs of $\cg$ such that $\Delta(N), \Delta(D_m) \leq C$ for some constant $C$. Let $R$ be a uniformly random perfect matching of $\cg - N$. Define $\bd X = (X_m)_{m=1}^\ell$, $X_m = |E(R \cap D_m)|$, $\bd Y = (Y_m)_{m=1}^\ell$ where $Y_m$ follows independently to $\po(|E(D_m)|/2n)$. Then,
    \[
    \lim_{n \to \infty} d_{TV}(\bd X, \bd Y) = 0
    \]
\end{corollary}
\begin{proof}
    Let $R^*$ be a uniformly random perfect matching from $\cg$. Consider the $\ell+1$-dimensional random vectors $(\bd X^*, X^*) = (X^*_1, ..., X^*_\ell, X^*), (\bd Y, Y^*) = (Y_1, ..., Y_\ell, Y^*)$ where $\bd X^* = (|E(R^* \cap D_m)|)_{m=1}^\ell$, and $X^* = |E(R^* \cap  N)|$, $Y^* \sim \po(|E(N)|/2n)$, $Y^*$ independent to $\bd Y$. We want to show the distribution $d_{TV}(\bd X^*| X^* = 0, \bd Y) \to 0$. We have
    \begin{align*}
    d_{TV}(\bd X^*|X^* = 0, \bd Y) & = \sum_{\bd k}|\pr(\bd X^* = \bd k|X^* = 0) - \pr(\bd Y = \bd k)| \\
    &= \sum_{\bd k}\left|\frac{\pr(\bd X^* = \bd k, X^* = 0)}{\pr(X^* = 0)} - \frac{\pr(\bd Y = \bd k)\pr(Y^* = 0)}{\pr(Y^* = 0)}\right| \\
    &=\sum_{\bd k}\left|\frac{\pr(\bd X^* = \bd k, X^* = 0)}{\pr(X^* = 0)} - \frac{\pr(\bd Y = \bd k, Y^* = 0)}{\pr(Y^* = 0)}\right| \qquad \text{by independence}\\
    &\leq \frac{1}{\pr(X^* = 0)}\left(\sum_{\bd k}|\pr(\bd X^* = \bd k, X^* = 0) - \pr(\bd Y = \bd k, Y^* = 0)|\right)  \\
    &+ \left|\frac{1}{\pr(X^* = 0)} - \frac{1}{\pr(Y^* = 0)}\right|\left(\sum_{\bd k}|\pr(\bd Y = \bd k, Y^* = 0)|\right) \\
    & \leq \frac{1}{\pr(X^* = 0)}d_{TV}((\bd X^*, X^*), (\bd Y, Y^*)) + \left|\frac{\pr(Y^* = 0)}{\pr(X^* = 0)}-1\right|
    \end{align*}
    We can use Theorem \ref{MainThm3.1} to show $d_{TV}((\bd X^*, X^*), (\bd Y, Y^*)) \to 0$ and use Corollary \ref{Cor3.2} to show $\lim_{n \to \infty} d_{TV}(\bd X^*|X^* = 0, \bd Y) = 0$.
\end{proof}

\section{Case of $K_{r \times (r-1)n}$}
\label{Sec4}
In this section, we assume $r$ is fixed, $\lambda= \binom r 2$. We will develop analogous results for Theorem \ref{MainThm3.1} in Section \ref{Sec3}, where $G = K_{r \times (r-1)n}$. The following theorem is a restatement of Theorem \ref{MainThm2} and is the main result of this section.
\begin{theorem}
    \label{MainThm4.1}
    Let $K_{r \times (r-1)n}$ be the balanced complete $r$-partite graph, let $(D_m)_{m=1}^\ell$ be a collection of disjoint subgraphs of $K_{r \times (r-1)n}$ with maximum degree $\Delta(D_m) \leq C$ for all $m$, where $C$ is independent of $n$. Let $R$ be a uniformly random \textbf{balanced} perfect matchings of $K_{r \times (r-1)n}$. Define random vector $\bd X = (|E(R \cap D_m)|)_{m=1}^\ell$, $\bd Y = (Y_m)_{m=1}^\ell$ where $Y_m \sim \po(|E(D_m)|/((r-1)^2n))$ independently. Then,
    \[
    \lim_{n \to \infty} d_{TV}(\bd X, \bd Y) = 0
    \]
\end{theorem}

We will first develop some notation for this section. Let $\bd x = (x_{m,i,j})_{m \in [\ell], i< j \in [r]} \in \mathbb N^{\ell \times \binom r 2}$. We define by convention that $x_{m,i,j} = x_{m,j,i}$ if $i>j$. For fixed $m$, we define $\bd x_m = (x_{m,i,j})_{i<j \in [r]} \in \mathbb N^{\binom r 2}$, for fixed $i,j$, we define $\bd x_{i,j} = (x_{m,i,j})_{m \in [\ell]} \in \mathbb N^\ell$, and for fixed $i$, we define $\bd x_i = (x_{m,i,j})_{m \in [\ell], j \neq i \in [r]} \in \mathbb N^{\ell \times (r-1)}$. We define $\vn \cdot$ as the sum of the coordinates as usual. 

We define $\psa: \mathbb N^{\ell \times \binom r 2} \to \mathbb N^{\ell}, \psb: \mathbb N^{\ell \times \binom r 2} \to \mathbb N^{r}, \psc: \mathbb N^{\ell \times \binom r 2} \to \mathbb N^{\binom r 2}$ as follows.
\begin{align*}
   \psa(\bd x) = (|\mathbf{x}_m|_1)_{m=1}^\ell \qquad \psb(\bd x) = (|\mathbf{x}_i|_1)_{i=1}^r \qquad \psc(\bd x) = (|\mathbf{x}_{i,j}|_1)_{i<j \in [r]}
\end{align*}
Let $V(K_{r \times (r-1)n}) = V_1 \uplus, ..., \uplus V_r$ be the vertex partition. We define the vector 
\[
\bd d = (d_{m,i<j})_{m \in [\ell], i<j \in [r]} \in \mathbb N^{\ell \times \binom r 2}
\] where $d_{m,i,j}$ is the number of edges of $D_m$ between $V_i$ and $V_j$. Let $\bd d_m = (d_{m,i,j})_{i<j \in [r]} \in \mathbb N^{\binom r 2}$. Then, $\vn{d_m} = |E(D_m)|$. The following definition will be analogous to Definition \ref{xmatching1}.
\begin{definition}
    Let $\bd x \in \mathbb N^{\ell \times \binom r 2}$. We define a $\bd x$-matching of $\cbg$ as a $\vn x$-matching where the number of edges in $D_m$ between $V_i$ and $V_j$ is $x_{m,i,j}$. We define $\mu_{\bd x}$ as the number of $\bd x$-matchings on $\cbg$.
\end{definition}
The next result is due to Johnston, Kayll and Palmer \cite{JKP}.
\begin{lemma}
    \label{bpmcount}
    The number of balanced perfect matching of $K_{r \times (r-1)n}$ is
    \[
    \bpm(\cbg) = \left(\frac{((r-1)n)!}{n!^{r-1}}\right)^r\left(n!\right)^{\binom r 2}
    \]
\end{lemma}
The following lemma is analogous to Lemma \ref{boundmu3}.
\begin{lemma}
    Let $\bd X$ be defined as Theorem \ref{MainThm4.1}. The probability generating function of $\bd X$ is given by
    \[
    G_{\bd X}(\bd s) = \sum_{\substack{\bd x \in \mathbb N^{\ell \times \binom r 2} \\ \forall i\neq j, \vn{x_{i,j}} \leq n}} \mu_{\bd x} \frac{(n \bd 1)_{(\psc(\bd x))}}{((r-1)n\bd 1)_{(\psb(\bd x))}} (\bd s - \bd 1)^{\psa(\bd x)}
    \]
\end{lemma}
\begin{proof}
    Let $\Omega$ be the set of all balanced perfect matching with a uniform probability measure, and let $\bd I = (I_{m,i,j})_{m \in [\ell], i<j \in [r]}$ where $I_{m,i,j}$ is the set of edges of $D_m$ between $V_i$ and $V_j$. For $\bd S \sqsubset \bd I$, let $\hat{\bd S}= \bigcup_{m \in [\ell],i<j \in [r]} S_{m,i,j}$ be the natural identification of $\bd S$ to a collection of edges in $\cbg$. Suppose $\hat{\bd S}$ is a submatching of a balanced perfect matching of $\cbg$. Then let $\bd x \in \mathbb N^{\ell \times \binom r 2}$ be such that $\hat{\bd S}$ is a $\bd x$-matching. We must have $\vn{x_{i,j}} \leq n$ for all $i\neq j \in [r]$. We want to count the number of ways to extend $\hat{\bd S}$ to a balanced perfect matching. Denote $V_i' = V_i - V(\hat{\bd S})$, and $|V_i'| = ((r-1)n - \vn{x_i})$. We want to partition $V_i'$ into $\{V_{i,j}\}_{j \neq i \in [r]}$  with $|V_{i,j}| = n - \vn{x_{i,j}}$. There are
    \[
    \frac{((r-1)n - \vn{x_i})!}{\prod_{j \neq i}(n - \vn{x_{i,j}})!}
    \]
    partitions. Then, we match the vertices in $V_{i,j}$ with the vertices in $V_{j,i}$ in $(n - \vn{x_{i,j}})!$ ways. Therefore, the probability equals
    \begin{align*}
     \pr(\hat{\bd S} \subset R) &= \pr(A_{\bd S}) \\
     &= \prod_i \frac{((r-1)n - \vn{x_i})!}{\prod_{j \neq i}(n - \vn{x_{i,j}})!}\prod_{i<j} (n - \vn{x_{i,j}})!\bigg/ \left(\frac{((r-1)n)!}{n!^{r-1}}\right)^r\left(n!\right)^{\binom r 2} \\
     &=\frac{(n \bd 1)_{(\psc(\bd x))}}{((r-1)n\bd 1)_{(\psb(\bd x))}}
    \end{align*}
    If $\hat{\bd S}$ is not a submatching of a balanced perfect matching, then $\pr(A_{\bd S}) = 0$. Since there are $\mu_{\bd x}$ $\bd x$-matchings, by Corollary \ref{pieCor}, we know
    \[
    G_{\bd X}(\bd s) \sum_{\substack{\bd x \in \mathbb N^{\ell \times \binom r 2} \\ \forall i\neq j, \vn{x_{i,j}} \leq n}} \mu_{\bd x} \frac{(n \bd 1)_{(\psc(\bd x))}}{((r-1)n\bd 1)_{(\psb(\bd x))}} (\bd s - \bd 1)^{\psa(\bd x)}
    \]
\end{proof}
The next lemma is analogous to Lemma \ref{boundmu3} and estimates the quantity $\mu_{\bd x}$.
\begin{lemma}
    Let $\mu_{\bd x}, C$ be defined as Theorem \ref{MainThm4.1}. Then,
    \[
    \mu_{\bd x} = \frac{(\bd d - O(\vn x)\bd 1)^{\bd x}}{\bd x!}
    \]
\end{lemma}
\begin{proof}
    We want to count the number of ways to choose the $\bd x$-matching. We assign an order to $(m,i,j): m \in [\ell], i<j \in [r]$ and choose the edges of $D_m$ between $V_i$ and $V_j$ following that order. Given $m_0, i_0, j_0$, suppose for all $m,i,j$ that precedes $m_0, i_0, j_0$ in that order, we have chosen the edges in the $\bd x$-matching that intersect $D_m$ and are between $V_i$ and $V_j$. Then, we want to select $x_{m_0,i_0,j_0}$ non-incident edges among the $d_{m_0,i_0,j_0}$ available edges. For $1 \leq k \leq x_{m_0,i_0,j_0}$, suppose we have already selected $k-1$ edges. Then, the total number of previously selected edges is no more than $\vn x$, and hence they are incident to no more than $2C \vn x = O(\vn x)$ edges because the maximum degree of $D_m$ is bounded by $C$. Therefore, the number of ways to select edges in $D_{m_0}$ between $V_{i_0}$ and $V_{j_0}$ is
    \[
    \frac{(d_{m_0,i_0,j_0}-O(\vn x))}{x_{m_0,i_0,j_0}!}
    \]
    where we divide by $x_{m_0,i_0,j_0}!$ because we do not distinguish orders among the $x_{m_0,i_0,j_0}$ edges.
    Finally, we take the product of all $m_0,i_0,j_0$ and by product rule, we obtained the expression by replacing the index with $m,i,j$.
\end{proof}
The next lemma will be used to bound the coefficients of the generating function.
\begin{lemma}
\label{boundcoef}
Given $n \geq \vn{x_{i,j}}$ for each $i \neq j \in [r]$, we have
\[
\frac{(n \bd 1)_{(\psc(\bd x))}}{((r-1)n\bd 1)_{(\psb(\bd x))}} \leq \frac{2^{2 \vn x}}{n^{\vn x}}
\]
\end{lemma}
\begin{proof}
    The lemma is equivalent to
    \[
    \frac{\prod_{i<j} n^{\vn{x_{i,j}}}n_{(\vn{ x_{i,j}})}}{\prod_i ((r-1)n)_{(\vn{x_{i}})}} \leq 2^{2 \vn x}
    \]
    Squaring on both sides, it suffices to show
    \[
    \prod_i \frac{\prod_{j \neq i} n^{\vn{x_{i,j}}}n_{(\vn{x_{i,j}})}}{((r-1)n)_{(\vn{x_{i}})}^2} \leq 2^{4 \vn x}
    \]
    Since for each $a \in \mathbb N, n^an_{(a)} \leq (2n)_{(2a)}$, it suffices to show that for each $i$,
    \[
    \frac{\prod_{j \neq i} (2n)_{(2\vn{x_{i,j}})}}{((r-1)n)_{(\vn{x_{i}})}^2} \leq 2^{2 \vn{x_i}}
    \]
    Since
    \[
    \frac{(2n)_{(2a)}}{n_{(a)}^2} \leq 2^{2a}
    \]
    we know
    \[
    \frac{\prod_{j \neq i} (2n)_{(2\vn{x_{i,j}})}}{((r-1)n)_{(\vn{x_{i}})}^2} \leq \prod_{j<i} \frac{(2n)_{(2\vn{x_{i,j}})}}{((r-j)n)^2_{(\vn {x_{i,j}})}} \prod_{j>i}\frac{(2n)_{(2\vn{x_{i,j}})}}{((r+1-j)n)^2_{(\vn {x_{i,j}})}} 
    \leq \prod_j \frac{(2n)_{(2\vn{x_{i,j}})}}{n_{(\vn{x_{i,j}})}^2} \leq 2^{2\vn{x_i}}
    \]
    and the proof follows.  
\end{proof}
Since $\bd Y = (Y_m)_{m=1}^\ell$ and $Y_m \sim \po(\vn{d_m}/(r-1)^2n)$, we can express the generating function of $\bd Y$ as
\begin{align*}
    G_{\bd Y}(\bd s) &= \exp\left(\frac{1}{(r-1)^2n}\psa(\bd d)\cdot (\bd s - \bd 1)\right) \\
    &= \sum_{\bd y \in \mathbb N^{\ell}} \frac{(1/((r-1)^2n))\psa(\bd d)^{\bd y}}{\bd y!}(\bd s - \bd 1)^{\bd y} \\
    &= \sum_{\bd x \in \mathbb N^{\ell \times \binom r 2}} \frac{(1/((r-1)^2n))\bd d)^{\bd x}}{\bd x!}(\bd s - \bd 1)^{\psa(\bd x)}
\end{align*}
where in the last step we used the multinomial theorem.

Set
\begin{align*}
    \alpha_{\bd x} &= \begin{cases}
        \mu_{\bd x} \frac{(n \bd 1)_{(\psc(\bd x))}}{((r-1)n\bd 1)_{(\psb(\bd x))}} & \text{for all} \ i\neq j \in [r], \vn{x_{i,j}} \leq n \\
        0 & \text{otherwise}
    \end{cases} \\
    \beta_{\bd x} &= \frac{(1/(r-1)^2n)\bd d)^{\bd x}}{\bd x!}
\end{align*}
By the same argument in Section \ref{Sec3}, we need to find $\gamma_{\bd x}$, independent of $n$, such that $0 \leq \alpha_{\bd x}(n), \beta_{\bd x}(n) \leq \gamma_{\bd x}$ and $\sum_{\bd x} \gamma_{\bd x}2^{\vn x} < \infty$. We also need to show 
\[
\lim_{n \to \infty} |\alpha_{\bd x}-\beta_{\bd x}| = 0
\]
holding $\bd x$ fixed. Then, by Tannery's Theorem, we have
\[
\lim_{n \to \infty} \dtv \leq \lim_{n \to \infty} \sum_{\bd x}|\alpha_{\bd x}-\beta_{\bd x}|2^{\bd x} = \sum_{\bd x} \lim_{n \to \infty} |\alpha_{\bd x}-\beta_{\bd x}|2^{\bd x} = 0
\]
\begin{proof}[Proof of Theorem \ref{MainThm4.1}]
We first want to find $\gamma_{\bd x}$ such that $\sum_{\bd x}\gamma_{\bd x}2^{\vn x} < \infty$.
By lemma \ref{boundcoef}
\[
\alpha_{\bd x}, \beta_{\bd x} \leq \frac{\bd d^{\bd x}} {\bd x!}\frac{2^{2 \vn x}}{n^{\vn x}} = \frac{((4/n)\bd d)^{\bd x}}{\bd x!}
\]
Since $d_{m,i,j}/n \leq (r-1)C$, we know
\[
\frac{((4/n)\bd d)^{\bd x}}{\bd x!} \leq \frac{(4(r-1)C)^{\vn x}}{\bd x!}
\]
and
\begin{align*}
    \sum_{\bd x \in \mathbb N^{\ell \times \binom r 2}} \frac{(4(r-1)C)^{\vn x}}{\bd x!} = \exp\left(4(r-1)C \binom r 2\right) < \infty
\end{align*}
Therefore, $\gamma_{\bd x} = (4(r-1)C)^{\vn x}/\bd x!$ will suffice. Then, we want to show $\lim_{n \to \infty}|\alpha_{\bd x}-\beta_{\bd x}| = 0$ for all fixed $\bd x$. By the triangle inequality, we can bound $|\alpha_{\bd x}-\beta_{\bd x}|$ as
\[
|\alpha_{\bd x}-\beta_{\bd x}| = \left|\frac{(\bd d - O(\vn x)\bd 1)^{\bd x}}{\bd x!} - \frac{\bd d^{\bd x}}{\bd x!}\right| \frac{(n \bd 1)_{(\psc(\bd x))}}{((r-1)n\bd 1)_{(\psb(\bd x))}} + \frac{\bd d^{\bd x}}{\bd x!}\left|\frac{(n \bd 1)_{(\psc(\bd x))}}{((r-1)n\bd 1)_{(\psb(\bd x))}} - \frac{1}{((r-1)^2n)^{\vn x}}\right|
\]
Applying Lemma \ref{boundprod} and Lemma \ref{boundcoef}, we know
\[
\left|\frac{(\bd d - O(\vn x)\bd 1)^{\bd x}}{\bd x!} - \frac{\bd d^{\bd x}}{\bd x!}\right| \frac{(n \bd 1)_{(\psc(\bd x))}}{((r-1)n\bd 1)_{(\psb(\bd x))}} \leq \ell \binom r 2 O(\vn x)((r-1)Cn)^{\vn x - 1}\frac{2^{2 \vn x}}{n^{\vn x}} = o(1)
\]
and
\[
\frac{\bd d^{\bd x}}{\bd x!}\left|\frac{(n \bd 1)_{(\psc(\bd x))}}{((r-1)n\bd 1)_{(\psb(\bd x))}} - \frac{1}{((r-1)^2n)^{\vn x}}\right| \leq ((r-1)Cn)^{\vn x}o\left(n^{-\vn x}\right) = o(1)
\]
holding $\bd x$ fixed. Therefore, $|\alpha_{\bd x}-\beta_{\bd x}| \to 0$ and the proof follows.
\end{proof}

\section{Graph Decomposition and case of non-disjont $D_m$}
\label{Sec5}
In this Section, we generalize $D_m$ to a collection of not necessarily disjoint subgraphs of $G$, where $G = \cg - N$ or $\cbg$. The idea is to decompose the $\ell$ subgraphs into $2^\ell -1$ disjoint pieces and apply the previous results to these pieces to get a $2^\ell-1$ independent Poisson joint distribution, and then appropriately sum the random variables with the corresponding distribution. We shall formalize this notion of graph decomposition using the following definition.
\begin{definition}
    Let $\dseq$ be a collection of subgraphs of $G$. Let $\mathcal P^*([\ell])$ be the collection of non-empty subsets of $[\ell]$. For each $S \in \mathcal P^*([\ell])$, define
    \[
    D_S = \bigcap_{m \in S} D_m \bigg / \bigcup_{m \in S^c} D_m^c
    \]
    Then, $(D_S)_{S \in \mathcal P^*([\ell])}$ is a collection of disjoint subgraphs.
\end{definition}
If we define the random variables $\bar X_S = |E(R \cap D_S)|$, then by applying our previous results, $\bar X_S$ approaches an independent Poisson joint distribution. Since $X_m = \sum_{S: m \in S} \bar X_S$ and the sum of a Poisson distribution is also Poisson, the random vector $(X_m)_{m=1}^\ell$ should approach a joint Poisson distribution. The next lemma characterizes this phenomenon.
\begin{lemma}
    \label{graphdecom}
     Consider four random vectors $\bd X = (X_m)_{m=1}^\ell, \bd Y = (Y_m)_{m=1}^\ell, \bar {\bd X} = (\bar X_S)_{S \in \mathcal P^*([\ell])}$, $(\bar Y_S)_{S \in \mathcal P^*([\ell])}$ such that $X_m = \sum_{S: m \in S} \bar X_S, Y_m = \sum_{S: m \in S} \bar Y_S$ for each $m$. Then,
     \[
     d_{TV}(\bd X, \bd Y) \leq d_{TV}(\bar {\bd X}, \bar{\bd Y})
     \]
\end{lemma}
\begin{proof}
    For vector $\bar {\bd k} = (\bar k_S)_{S \in \mathcal P^*([\ell])}$, let $\bd k = (k_m)_{m=1}^\ell$ be defined as $k_m = \sum_{S: m \in S} \bar k_S$. We define a function $T : \mathbb N^{\mathcal P^*([\ell])} \to \mathbb N^{\ell} $ such that $T(\bar {\bd k}) = \bd k$. The function $T$ is surjective, because if we let $\bar k_S = k_m$ if $S = \{m\}$ and $\bar k_S = 0$ otherwise, and let $\bar{\bd k} = (\bar k_S)_{S \in \mathcal P^*([\ell])}$, then $T(\bar{\bd k}) = (k_m)_{m=1}^\ell$. Then, we have
     \begin{align*}
         d_{TV}(\bd X, \bd Y) & =\sum_{\bd k \in \mathbb N^\ell}\left|\pr(\bd X = \bd k) - \pr(\bd Y = \bd k)\right| \\
         &= \sum_{\bd k \in \mathbb N^\ell}\left|\sum_{\bar{\bd k} \in T^{-1}(\bd k)} \pr(\bar{\bd X} = \bar{\bd k})-\pr(\bar{\bd X} = \bar{\bd k})\right| \\
         &\leq \sum_{\bd k \in \mathbb N^\ell}\sum_{\bar{\bd k} \in T^{-1}(\bd k)}\left|\pr(\bar{\bd X} = \bar{\bd k})-\pr(\bar{\bd X} = \bar{\bd k})\right| \\
         &= \sum_{\bar {\bd k} \in \mathbb N^{\mathcal P^*([\ell])}}\left|\pr(\bar{\bd X} = \bar{\bd k})-\pr(\bar{\bd X} = \bar{\bd k})\right| \\
         &= d_{TV}(\bar{\bd X}, \bar{\bd Y})
     \end{align*}
\end{proof}
The next corollary generalizes Theorem \ref{MainThm1}, and the proof follows naturally from this lemma and Theorem \ref{MainThm1}.
\begin{corollary}
    Let $N$ be a subgraph of $\cg$, and let $(D_m)_{m=1}^\ell$ be a collection of subgraphs of $\cg - N$ such that $\Delta(D_m), \Delta(N) \leq C$ for some constant $C$. Let $(D_S)_{S \in \mathcal P^*([\ell])}$ be the graph decomposition of $(D_m)_{m=1}^\ell$. Let $R$ be a uniformly chosen random perfect matching of $\cg - N$ and define $\bd X = (X_m)_{m=1}^\ell, X_m = |E(R \cap D_m)|$. For each $S \in \mathcal P^*([\ell])$, let $Y_S$ follow independently to $\po(|E(D_S)|/2n)$ and let $Y_m = \sum_{S: m \in S}Y_S$, and $\bd Y = (Y_m)_{m=1}^\ell$. Then
    \[
    \lim_{n \to \infty} d_{TV}(\bd X, \bd Y) = 0
    \]
\end{corollary}
\begin{proof}
  We know $(D_S)_{S \in \mathcal P^*([\ell])}$ is a collection of disjoint subgraphs of $\cg - N$ where $\Delta(D_S) \leq C$ for all $S$. Define $\bar{\bd X} = (|E(R \cap D_S)|)_{S \in \mathcal P^*([\ell])}$, $\bar{\bd Y} = (Y_{S})_{S \in \mathcal P^*([\ell])}$. By Theorem \ref{MainThm1}, we know
\[
d_{TV}(\bar{\bd X},\bar{\bd Y}) \to 0
\]
Therefore, by Lemma \ref{graphdecom}, $\dtv \to 0$. 
\end{proof}
The next corollary generalizes Theorem \ref{MainThm2} and is proven by the exact same method. We merely replace $\cg - N$ by $\cbg$ and use Theorem \ref{MainThm2} instead of Theorem \ref{MainThm1}.
\begin{corollary}
    Let $\dseq$ be a collection of subgraphs of $\cbg$ such that $\Delta(D_m) \leq C$ for some constant $C$. Let $(D_S)_{S \in \mathcal P^*([\ell])}$ be the graph decomposition of $(D_m)_{m=1}^\ell$. Let $R$ be a uniformly chosen random perfect matching of $\cbg$ and define $\bd X = (|E(R \cap D_m|)_{m=1}^\ell$. For each $S \in \mathcal P^*([\ell])$, let $Y_S$ follow independently to $\po(|E(D_S)|/(r-1)^2n)$ and let $Y_m = \sum_{S: m \in S}(Y_S)$, and $\bd Y = (Y_m)_{m=1}^\ell$. Then
    \[
    \lim_{n \to \infty} d_{TV}(\bd X, \bd Y) = 0
    \]
\end{corollary}

\section{Potential Directions of Future Work}
\label{Sec6}

A promising direction for future research is to extend the analysis from the parent graph $G$ to a broader class of graphs, such as regular robust expander graphs (see \cite{GJ, RBexpander1, RBexpander2, RBexpander3} for details). The current method of counting via the Principle of Inclusion-Exclusion may not be well-suited for this generalization. Instead, adopting techniques similar to those employed in \cite{GJ} could provide a more effective approach.

Another potential direction for future research involves relaxing the constraints on the constants $ \ell $, $ C $, and $ r $, allowing them to grow, possibly slowly, as a function of $ n $. If $ \ell \to \infty $ as $ n \to \infty $, the current approach of applying Tannery’s theorem to interchange the limit and infinite summation becomes inapplicable, as the theorem requires a countable index set $ I $, whereas this scenario involves the uncountable set $ \mathbb{N}^{\mathbb{N}} $. Furthermore, the graph decomposition technique outlined in Section 5 may no longer be effective, as it increases the number of decomposed graphs $ D $ from $ \ell $ to $2^\ell - 1 $, where the latter grows significantly faster than the former.

A third direction for future research is to investigate the rate of convergence of $\dtv$. Current results establish only that $\dtv = o(1) $. We conjecture that, depending on the specific properties of the graph sequence $\dseq$, it may be possible to identify a constant $ K = K(\dseq) $ such that $\dtv = K n^{-1} + o(n^{-1})$.

\section{Acknowledgement}
The author thanks his advisor, Professor Karen L. Collins, in the Department of Mathematics and Computer Science at Wesleyan University, for advising this research process.

\end{document}